 \documentclass[12pt,reqno]{amsart}
 \usepackage[dvips]{epsfig}
 \usepackage{amsgen, amstext,amsbsy,amsopn, amsthm, amsfonts,amssymb,amscd,amsmat
 h,euscript,enumerate,url,verbatim,calc,xypic}
 \usepackage{latexsym}
 \usepackage{graphics}
 \usepackage{color}

 \newcommand{\n}{\mathfrak{n} }
 \newcommand{\m}{\mathfrak{m} }

 \newcommand{\p}{\mathcal{P}}

 \newcommand{\R}{\mathcal{R}}

    \newcommand{\coker}{\operatorname{coker}}

  \newcommand{\ann}{\operatorname{ann}}

 \newcommand{\dm}{\operatorname{dim}}
 \newcommand{\h}{\operatorname{ht}}
 
 \newcommand{\Supp}{\operatorname{Supp}}

 \newcommand{\depth}{\operatorname{depth}}

 \newcommand{\Tor}{\operatorname{Tor}}

  \newcommand{\lm}{\lambda}

\theoremstyle{plain}
 \newtheorem{theorem}{Theorem}[]

 \newtheorem{lemma}[theorem]{Lemma}

 \theoremstyle{definition}

 \theoremstyle{remark}

\title[ ] {On the Chern number of an ideal }
\author{Mousumi Mandal}
\author{J. K. Verma} 
\address{Department of Mathematics, IIT Bombay, Mumbai 400 076}
\email{mousumi@math.iitb.ac.in}
\email{jkv@math.iitb.ac.in}
                    

\begin{document}
\maketitle
 
\begin{abstract} We settle the negativity conjecture of Vasconcelos for the Chern number of an ideal in certain unmixed quotients of regular local rings by explicit calculation of the Hilbert polynomials  of
all ideals generated by systems of parameters.

\end{abstract}

\section*{Introduction}
Let $I$ be an $\m$-primary ideal in a local ring  $(R,\m)$ 
of dimension $d.$ Let $H(I,n)=\lm(R/I^n)$ denote the Hilbert function of $I$ where $\lm(M)$ denotes the length of an $R$-module $M.$  
The Hilbert function $H(I,n)$ is given by a polynomial $P(I,n)$ of degree $d$ for large $n.$ It is written in the form
$$P(I,n)=e_0(I){n+d-1 \choose d}-e_1(I){n+d-2 \choose d-1} + \cdots
+(-1)^d e_d(I).$$

If $I$ is generated by a system of parameters, then $R$ is Cohen-Macaulay if and only if $e_0(I)=\lm(R/I).$ Recently it has been observed by Vasconcelos \cite{vas} that the signature of the 
coefficient $e_1(I),$ called the {\em Chern number} of $I,$
can be used to characterize Cohen-Macaulay property of $R$ for 
large classes of rings. In the Yokohama Conference in 2008 Vasconcelos  proposed   the following: \\

\noindent
{\bf  Negativity Conjecture:}  Let $(R,\m)$ be  an unmixed, equidimensional local ring which  is a homomorphic image
of a Cohen-Macaulay local ring. Then $R$ is not Cohen-Macaulay
if and only if for any ideal $J$ generated by a system of parameters, $e_1(J) < 0.$

Ghezzi, Hong and Vasconcelos \cite{ghv} settled the conjecture for  (1) Noetherian local domains of dimension $d \geq  2,$ which are  homomorphic images of  Cohen-Macaulay local rings and  for (2) Universally catenary integral domains containing 
a field.  The Negativity Conjecture has been resolved
for all unmixed local rings by S. Goto recently. 
    
In this paper  we settle the Negativity Conjecture for certain unmixed quotients of regular local rings, by explicitly finding the  Hilbert polynomial of  all parameter ideals. 

{\bf Acknowledgements:} We wish to thank Shiro Goto, A. Simis and V. Srinivas for several useful conversations.

\section{\bf The Hilbert polynomial of parameter ideals  in certain  quotients  of  regular local rings}
L. Ghezzi, J. Hong and W. Vasconcelos \cite{ghv} calculated
the Chern number of any parameter  ideal in certain quotients 
of regular local rings of dimension four. We recall their  result first: 

\begin{theorem} Let $(S,\m)$ be a four dimensional
regular local ring with $S/\m$ infinite. Let $P_1, P_2, \ldots, P_r$ be a family of height two prime ideals of $S$ so that for $i \neq j, P_i+P_j$ is $\m$-primary. Put $R=S/\cap_{i=1}^{r} P_i.$ Let
$J$ be an $\m$-primary  parameter ideal  of $R.$
Let $L= \oplus_{i=1}^r S/P_i/R.$ If $J \subseteq \ann L$ then $e_1(J)=-\lm(L)$ and $e_2(J)=0.$
\end{theorem}

In this section we find the Hilbert polynomials  of parameter ideals   in $R$ when $S$ is any regular local ring.

\begin{lemma}
 Let $(S,\n)$ be an $r$-dimensional regular local ring. Let $I$ be a height $h$ Cohen-Macaulay ideal of $S$.  Suppose  $a_1,\ldots ,a_d\in S$  such that  $(a_1+I,\ldots, a_d+I)$ is a system of parameters in $S/I$. Then $a_1,\ldots , a_d$ is a regular sequence in $S$.  
\end{lemma}
\begin{proof}
 Let $J=(a_1,\ldots ,a_d)$. Then $\lm(S/I\otimes_S S/J)<\infty $.
Hence by Serre's theorem \cite[Theorem 3, Chapter 5]{ser} 
 $$\dm S/I +\dm S/J\leq \dm S=r.$$
 As $S$ is regular, it is catenary. Thus $r-\h I+r-\h J\leq r$. Therefore $d\leq \h J\leq d$. Hence $\h J=d$ and consequently $a_1,\ldots ,a_d$ is an $S$-regular sequence.
\end{proof}
\begin{lemma}\label{lt}
 Let $J,I$ and $S$ be as above. Then for all $j,n\geq 1$ $$\Tor_j^S(S/J^n,S/I)=0.$$
\end{lemma}
\begin{proof}
 We will apply induction on $n$. Let $n=1$. As $J$ is a complete intersection, the Koszul complex $K(\underbar{a})$ of the sequence
$\underbar{a}=a_1, a_2,\ldots, a_d$  
 $$K(\underbar{a}): 0\longrightarrow S\longrightarrow S^d\longrightarrow S^{d\choose 2}\longrightarrow \ldots \longrightarrow S^d\longrightarrow S\longrightarrow S/J\longrightarrow 0$$
gives  a free resolution of $S/J$. Tensoring the above complex  with $R:=S/I$ we get
 $$K(\underbar{a},S/I):0\longrightarrow R\longrightarrow R^d\longrightarrow R^{d\choose 2}\longrightarrow \ldots \longrightarrow R^d\longrightarrow R\longrightarrow R/K\longrightarrow 0$$ which is the Koszul complex of $JR=K$. As $K$ is  generated by an $R$-regular sequence, the above is a free resolution of $R/K$. Hence $\Tor_j^S(S/J,S/I)=0$ for all $j\geq 1$. Since $J$ is generated by a regular sequence, $J^n/J^{n+1}$ is a free $S/J$-module. Consider the exact sequence 
 $$0\longrightarrow J^n/J^{n+1}\longrightarrow S/J^{n+1}\longrightarrow S/J^n\longrightarrow 0.$$
 This gives rise to the long exact sequence
 $$\ldots \longrightarrow \Tor_j^S(J^n/J^{n+1},S/I)\longrightarrow \Tor_j^S(S/J^{n+1},S/I)\longrightarrow \Tor_j^S(S/J^{n},S/I)\longrightarrow \ldots $$
 By induction on $n$, it follows that $\Tor_j^S(S/J^{n+1},S/I)=0$ for all $j\geq 1$.
\end{proof}
\begin{lemma}
 Let $J=(a_1,\ldots,a_d)$ be a complete intersection of height $d$ in a regular local ring $(R,\m)$. Let $L$ be an $R$-module of finite length. Then $\R(J)\otimes_R L$ is a finite $\R(J)$-module of dimension $d$ and 
 $$\Supp(\R(J)\otimes_R L)=V(\m\R(J)).$$ 
\end{lemma}
\begin{proof}
 Suppose $\p \in \Supp(\R(J)\otimes_R L)$ and $p=\p\cap R$. Then $(\R(J)\otimes_R L)_{\p}=(\R(J_p)\otimes_{R_p}L_p)_\p\not= 0$. Hence $L_p\not= 0$. As $L$ has finite length, $p=\m$. Hence $\m\R(J)\subseteq \p$. Since $\R(J)/\m\R(J)\simeq R/\m[T_1,T_2,\ldots, T_d]$ is a polynomial ring, $\m\R(J)$ is a prime ideal of $\R(J).$

 We will prove the other inclusion by induction on $\lm(L)$. Let $\lm(L)=1$. Then $L\cong R/\m$. Therefore $\R(J)\otimes_R L=\R(J)\otimes_R R/\m=\R(J)/\m \R(J).$ Hence $\m\R(J)\in \Supp(\R(J)\otimes_R L)$. Now assume that $\lm(L)>1$. Then we have the following exact sequence of $R$-modules
 \begin{equation}\label{eq2}
 0\longrightarrow R/\m\xrightarrow{i}L\longrightarrow C\longrightarrow 0
 \end{equation}
 where $C=\coker i$. Tensoring (\ref{eq2}) with $\R(J)$ we get the exact sequence
 \begin{equation}
   \R(J)\otimes_R R/\m\longrightarrow\R(J)\otimes_R L\longrightarrow \R(J)\otimes_R C\longrightarrow 0.
 \end{equation}
 Localize the above sequence at $\m\R(J)$  to get  $\m\R(J)\in\Supp(\R(J)\otimes_R L)$ using induction.
 \end{proof}

\begin{lemma}\label{ko}
Let $S$ be a local ring, $a_1,\ldots ,a_d$ be a regular sequence and $J=(a_1,\ldots ,a_d).$ Let $L$ be an $S$-module of finite length. 
If $J\subseteq \ann L$ then 
$$\lm(\Tor_1(L,S/J^n))={n+d-1\choose d-1}\lm(L).$$
\end{lemma}
\begin{proof}
 By \cite[Example 10]{kod}  for any $n>0$, $J^n$ is generated by the maximal minors of the $n\times (n+d-1)$ matrix $A$ where 
 \[
  A=\left( \begin{array}{ccccccccc}
            a_1 & a_2 & a_3 &\cdots & a_d & 0 & 0 & \cdots & 0\\
	    0   & a_1 & a_2 & \cdots & a_{d-1} & a_d & 0 & \cdots & 0\\
	    0   &  0  & a_1 & \cdots & a_{d-2} & a_{d-1} & a_d & \cdots & 0\\
	    0   &  0  &  0 & \cdots & a_1 & a_2 & a_3 & \cdots & a_d
           \end{array}
\right) 
 \]
By  Eagon-Northcott \cite[Theorem 2]{en}, the minimal free resolution of $S/J^n$ is given by
\begin{equation}\label{eg}
0\longrightarrow S^{\beta_{d}}\longrightarrow S^{\beta_{d-1}}\longrightarrow \cdots \longrightarrow S^{\beta_1}\longrightarrow S\longrightarrow S/J^n \longrightarrow 0
\end{equation}
where the Betti numbers of $S/J^n$ are given by $$\beta_i^S(S/J^n)={n+d-1\choose d-i}{n+i-2\choose i-1},~~~~~~~~~~~~1\leq i\leq k.$$
Taking tensor product of (\ref{eg}) with $L$ we get the following complex
$$0\longrightarrow L^{\beta_{d}}\longrightarrow L^{\beta_{d-1}}\longrightarrow \cdots \longrightarrow L^{\beta_1}\longrightarrow L\longrightarrow L/J^nL \longrightarrow 0.$$
Since $J\subseteq \ann L,$ the  maps in the above complex are zero. 
Hence 
$$\lm(\Tor_1(L,S/J^n))=\beta_1\lm(L)={n+d-1\choose d-1}\lm(L).$$
\end{proof}

\begin{theorem}
 Let $(S,\n)$ be a regular local ring of dimension $r$ and $I_1,\ldots ,I_g$ be Cohen-Macaulay ideals of height $h$ which satisfy the condition: $I_i+I_j$ is $\n$-primary for $i\not= j$. Let $R=S/I_1\cap \ldots \cap I_g$ and $d=\dm R\geq 2$. Let $a_1,\ldots ,a_d\in S$ such 
 that their images in $R$ form  a system of parameters. Let $J=(a_1,\ldots ,a_d),$  $L=[\oplus_{i=1}^gS/I_i]/R$ and  $K=JR$. 
Put $H_J(L,n)=\lm(J^n\otimes_R L)$ and  $P_J(L,n)$ be the corresponding Hilbert polynomial. Then

{\em (1)}\; $P_J(L,n)=-e_1(K){n+d-2\choose d-1}+e_2(K){n+d-3\choose d-2}- \cdots +(-1)^de_d(K)+\lm(L).$

{\em (2)}\;   If $ J\subseteq \ann L$  then 
  $$P(K,n)= e_0(K) {n+d-1 \choose d}+ \lm(L){n+d-2\choose d-1}+\lm(L){n+d-3\choose d-2}+\cdots +n\lm(L).$$

\end{theorem}
\begin{proof}

First we show that $\lm(L) < \infty$. Consider the exact sequence
 \begin{equation}\label{e1}
 0\longrightarrow \frac{S}{I_1\cap \cdots \cap I_g}\longrightarrow N=\frac{S}{I_1}\oplus \ldots \oplus \frac{S}{I_g}\longrightarrow L\longrightarrow 0.
 \end{equation}
 Let $P$ be a non-maximal prime ideal of $S$ not containing any $I_1,\ldots ,I_g$. Then $L_P=0$. If there is an $i$ such that $I_i\subseteq P$, then for $j\not= i$, $I_j\nsubseteq P$. Hence
\begin{equation*}
(S/I_1\cap \ldots \cap I_g)_P=(S/I_1\oplus \ldots \oplus S/I_g)_P=(S/I_i)_P.
\end{equation*}
Thus $L_P=0$. Hence $\Supp L=\{\n\}$. Thus $\lm(L)<\infty$.

By the depth lemma, 
$\depth R=1$. Thus $R$ is not Cohen-Macaulay. 
 Tensoring (\ref{e1}) with $S/J^n$  we get the exact sequence
 $$\longrightarrow \Tor_1^S(R,S/J^n)\longrightarrow \bigoplus_{i=1}^g\Tor_1^S(S/I_i,S/J^n)\longrightarrow\Tor_1^S(L,S/J^n)$$
 $$\longrightarrow R/K^n\longrightarrow \bigoplus_{i=1}^g S/(I_i,J^n)\longrightarrow L/J^nL\longrightarrow 0.$$
 By Lemma \ref{lt}, $\Tor_1^S(S/I_i,S/J^n)=0$ for all $i,n$. For large $n$, $J^nL=0$ as $\lm(L) < \infty$. For large $n,$
 $$\lm(R/K^n)=e_0(K){n+d-1\choose d}-e_1(K){n+d-2\choose d-1}+\cdots +(-1)^de_d(K).$$
 By (\ref{e1}) and additivity of $e_0(J,\_)$ we get $e_0(K)=\sum_{i=1}^ge_0(J,S/I_i)$.
 Hence $$\lm(\Tor_1^S(L,S/J^n))-\lm(L)=\sum_{i=1}^d
(-1)^i e_i(K){n+d-1-i\choose d-i}.$$ From the exact sequence
 $$0\longrightarrow J^n\longrightarrow S\longrightarrow S/J^n\longrightarrow 0$$
 we get
  $$\longrightarrow \Tor_1^S(J^n,L)\longrightarrow \Tor_1^S(S,L)\longrightarrow\Tor_1^S(S/J^n,L)$$
  $$\longrightarrow  J^n\otimes_S L\longrightarrow S\otimes L\longrightarrow L/J^nL\longrightarrow 0.$$
Hence for large $n$, 
\begin{equation}\label{eq3}  
\lm(\Tor_1^S(S/J^n,L))=\lm(J^n\otimes_SL)= \left[\sum_{i=1}^d (-1)^i e_i(K){n+d-1-i\choose d-i}\right]+\lm(L).
 \end{equation}
Since $\dm(\R(J)\otimes_RL)=d$ and $d\geq 2$, $e_1(K)<0$.
 If $J\subseteq\ann L$ then by Lemma \ref{ko} $$\lm(\Tor_1(L,S/J^n))={n+d-1\choose d-1}\lm(L).$$
Substituting in (\ref{eq3}) we get
$${n+d-1\choose d-1}\lm(L)-\lm(L)=-e_1(K){n+d-2\choose d-1}+e_2(K){n+d-3\choose d-2}- \cdots +(-1)^de_d(K).$$
Using the  equation,
$$ {n+d-1 \choose d-1}= 1+ \displaystyle{\sum_{i=1}^{d-1}{n+d-i-1\choose d-i}} $$
we obtain $e_i(K)=(-1)^i\lm(L)$ for $i=1,2,\ldots, d-1$ and $e_d(K)=0$.

\end{proof}


\begin{thebibliography}{AAAA}

\bibitem {bh} W. Bruns and J. Herzog, Cohen-Macaulay Rings, Revised Edition, 
Cambridge University Press, 1998.

\bibitem {en} J. A. Eagon and D. G. Northcott,
{\em Ideals defined by  matrices and a certain complexe associated with them }, Proc. Royal Soc. Ser. A {\bf 269} (1962) 188-204.

\bibitem {ghv} L. Ghezzi, J. Hong and W. Vasconcelos, 
{\em The signature of the Chern coefficients of local rings}, arXiv:0807.2686v1, 17 July 2008.

\bibitem {gn} S. Goto and K. Nishida,
{\em Hilbert coefficients and Buchsbaumness of associated graded rings},
J. Pure Appl. Algebra {\bf 181} (2003) 61--74.


\bibitem {kod} V. Kodiyalam, 
{\em Homological invariants of powers of an ideal}, 
Prof. Amer. Math. Soc. {\bf 118} (1993) 757-764.

\bibitem {ser} J.-P. Serre, Local Algebra, Springer-Verlag (Berlin), 2000.

\bibitem  {vas} W. Vasconcelos, {\em The Chern coefficients of local rings}, Michigan Math.  J. {\bf 57} (2008), 725-743. 
\end{thebibliography}
\end{document}